\documentclass[preprint,12pt]{elsarticle}

\usepackage{amssymb}
\usepackage{amscd,amsmath,amsfonts,graphicx}
\usepackage{amsthm}

\begin{document}
\baselineskip=22pt \centerline{\large \bf On Powers of the Catalan
Number Sequence} \vspace{1cm} \centerline{Gwo Dong Lin}
\centerline{Academia Sinica, Taiwan} \vspace{1cm} \noindent {\bf
Abstract.} The Catalan number sequence is one of the most famous
number sequences in combinatorics and is well studied in the
literature. In this paper we further investigate its fundamental
properties related to the moment problem and prove for the first
time that it is an infinitely divisible Stieltjes moment sequence in
the sense of S.-G.\,Tyan. Besides, any positive real power of the
sequence is still a Stieltjes determinate sequence. Some more cases
including (a) the central binomial coefficient sequence (related to
the Catalan sequence), (b) a double factorial number sequence and
(c) the generalized Catalan (or Fuss--Catalan) sequence are also
investigated. Finally, we pose two conjectures including the
determinacy equivalence between powers of nonnegative random
variables and powers of their moment sequences, which is supported
by some existing results.
\vspace{0.1cm}\\
\hrule
{\bf Short title: Powers of Catalan Numbers}\\
{\bf Postal address:} Gwo Dong Lin, Institute of Statistical
Science, Academia Sinica, Taipei 11529, Taiwan. (E-mail:
 gdlin@stat.sinica.edu.tw)\\ \bigskip\\
  \noindent\textsc{to appear in \emph{Discrete
Mathematics} (2018)} https://doi.org/10.1016/j.disc.2018.05.009\\
\newpage

\noindent{\bf 1. Introduction}
\newcommand{\bin}[2]{
   \left(
     \begin{array}{@{}c@{}}
         #1  \\  #2
     \end{array}
   \right)          }

The Catalan number sequence  was introduced by Catalan (1838) and
defined by $C_n$ = $2n\choose n$$/(n+1)=\frac{(2n)!}{n!(n+1)!},\
n=0, 1,2,\dots.$ It is one of the most famous and frequently
encountered  number sequences in combinatorics (Koshy 2009) and has
the generating function
$${C(x)=\frac{1-\sqrt{1-4x}}{2x}}=\frac{2}{1+\sqrt{1-4x}}=\sum_{n=0}^{\infty}C_nx^n,\ |x|<1/4$$
(Binet 1839). In the literature, there are more than $200$
interpretations of the sequence so far (Stanley 2015, Roman 2015).
It worths mentioning that An-Tu Ming (1692 -- 1763) also used this
sequence in the approximation theory, e.g., he derived the following
identity of the sine function in terms of $C_n$ (Larcombe 1999, Luo
1988, 2013):
$$\left(\sin{\frac{\alpha}{2}}\right)^2=\sum_{n=1}^{\infty}C_{n-1}\left(\frac{\sin
\alpha}{2}\right)^{2n},\ \ \alpha\in[0,\pi/2].$$

On the other hand, let us recall that a sequence
$\{m_n\}_{n=0}^{\infty}$ of real numbers is called a Stieltjes
moment sequence if there exists a {nonnegative measure $\mu$} on
${\mathbb R}_+\equiv[0,\infty)$ having it as its moment sequence:
\begin{eqnarray}
m_n=\int_{[0,\infty)}x^nd\mu(x),\ n=0,1,2,\ldots.
\end{eqnarray}
If, in addition, such a measure $\mu$ in (1) is unique, then
$\{m_n\}_{n=0}^{\infty}$ is called a  Stieltjes determinate (S-det)
sequence; otherwise, $\{m_n\}_{n=0}^{\infty}$ is called a Stieltjes
indeterminate (S-indet) sequence. Note that if $m_0=1,$ then (1) can
be rewritten as
\begin{eqnarray}
m_n=\int_0^{\infty}x^ndF(x),\ n=0,1,2,\ldots,
\end{eqnarray}
where $F$ is the distribution function of a random variable $X\ge
0,$ denoted $X\sim F.$ If $F$ in (2) is unique by moments (namely,
there is no other distribution having the same moment sequence
 as $F$), we also say that $F$ (or $X$) is
moment-determinate (M-det) on ${\mathbb R}_+$; otherwise, $F$ (or
$X$) is moment-indeterminate (M-indet) on ${\mathbb R}_+$.

 Since both the Hankel
matrices $\Delta_n=(C_{i+j})_{0\le i,j\le n}$ and
$\overline{\Delta}_n=(C_{i+j+1})_{0\le i,j\le n}$ of the Catalan
sequence have nonnegative determinants for all $n=0,1,2,\ldots,$
$\{C_n\}_{n=0}^{\infty}$ is a Stieltjes moment sequence by the
remarkable characterization result due to Stieltjes (1894/1895) (see
also Shohat and Tamarkin 1943, or Akhiezer 1965). More precisely, we
have
\begin{eqnarray*}
{\hbox{det}(\Delta_n)=  \hbox{det}(\overline{\Delta}_n)=1\
\forall n}
\end{eqnarray*} and this happens to be a characteristic property of the Catalan number sequence
(Stanley 1999, Chapter 6). Actually, the measure $\mu$ corresponding
to $\{C_n\}_{n=0}^{\infty}$  has a density $f_C,$ called Catalan
density and described below:
\begin{eqnarray}C_n=\int_0^{\infty}x^nd\mu(x)=\int_0^{\infty}x^nf_C(x)dx\equiv\int_0^{\infty}x^n
{\left(\frac{1}{2\pi}\sqrt{\frac{4-x}{x}}\right)}I_{(0,4]}(x)dx , \forall n,
\end{eqnarray}
where $I_A$ is the indicator function of the set $A$ (Stanley 2015,
Amdeberhan et al.\,2013). Namely, $\{C_n\}_{n=0}^{\infty}$ is
exactly the moment sequence of a bounded nonnegative random variable
$X_C\sim F_C$ with density $f_C,$ and $X_C$ is M-det on both
${\mathbb R}_+$ and ${\mathbb R}\equiv(-\infty,\infty),$ because its
moment generating function exists. (See, e.g., Cram\'er's and
Hardy's criteria in Lin 2017, Theorems 1 and 2.) Moreover, the
random variable $X_C/4$ has Beta distribution $B_{\alpha,\beta}$
with parameters $\alpha=1/2,\ \beta=3/2$ and moment sequence
$\{2^{-2n}C_n\}_{n=0}^{\infty}.$ The Catalan distribution $F_C$ is
also called the Marchenko-Pastur (or the free Poisson) distribution
in the literature (see, e.g., Banica et al.\,2011).

Recently, Liang et al.\,(2016) studied the conditions under which
the Catalan-like number sequences are Stieltjes moment sequences,
while Berg (2005, 2007) investigated the moment (in)determinacy
property of the powers of the factorial sequence
$\{n!\}_{n=0}^{\infty}$. Mimicking Berg's approach  (Lemma 3 below),
 we have some new findings on the topic. We first prove that {\it for
any real} $c>0,$ the $c$-th power of the Catalan sequence,
$\{C_n^c\}_{n=0}^{\infty},$ is a Stieltjes moment sequence, namely,
$\{C_n\}_{n=0}^{\infty}$ is an {\it infinitely divisible} Stieltjes
moment sequence (Tyan 1975), denoted $\{C_n\}_{n=0}^{\infty}\in
{\cal I}$, and then prove that each power sequence
$\{C_n^c\}_{n=0}^{\infty}$ is even S-det.

The five cases (i) the central binomial coefficient sequence,
$\{{2n\choose n}\}_{n=0}^{\infty},$ (ii) the {\it double factorial}
number sequence, $\{(2n-1)!!\}_{n=0}^{\infty},$  (iii) the
generalized Catalan (or Fuss--Catalan) sequence,  $\{{(k+1)n\choose
n}/(kn+1)\}_{n=0}^{\infty},$ (iv) $\{{(k+1)n\choose
n}\}_{n=0}^{\infty}$ and (v) $\{(kn)!\}_{n=0}^{\infty}$ are also
investigated.

The main results are stated in the next section, and their proofs
are given in Section 4. Section 3 provides some necessary lemmas.
Finally, in Section 5 we pose
 two conjectures including
the determinacy equivalence between powers of nonnegative random
variables and powers of their moment sequences, which is supported
by some existing results.
\bigskip\\
{\bf 2. Main results}

We start with a complete study in the simplest cases.

 \noindent{\bf Theorem 1.}
(a) The Catalan sequence
 is an infinitely divisible Stieltjes moment sequence, i.e., $\{C_n\}_{n=0}^{\infty}\in {\cal I}.$ \\
(b) For each real $c>0,$ the power sequence
$\{C_n^c\}_{n=0}^{\infty}$ is S-det.\\
(c) For each real $c>0,$ the measure $\mu_c$ (with $\mu_1=\mu$ in
(3)) corresponding to the moment sequence $\{C_n^c\}_{n=0}^{\infty}$
has the Mellin transform
$${\cal
M}_c(s)\equiv\int_0^{\infty}t^sd{\mu}_c(t)=\left(\frac{2^{2s}}{\sqrt{\pi}}\cdot\frac{\Gamma(s+1/2)}{\Gamma(s+2)}\right)^c
=\left(\frac{1}{s+1}\cdot \frac{\Gamma(2s+1)}{(\Gamma(s+1))^2}\right)^c,\
s\ge 0.$$

In the next theorem, we treat the central binomial coefficient
 sequence by a new approach, which is interesting in itself.
 See also Lemma 6 below, in which the equality in distribution happens to
 be a characteristic property of the density in (4), as shown in
 Remark 1.

 \noindent{\bf Theorem
2.} Let $B_0=1,\ B_n=$$2n\choose n$$=(n+1)C_n$,
$n=1,2,\ldots.$ Then we have\\
(a) the central binomial coefficient
 sequence $\{B_n\}_{n=0}^{\infty}\in {\cal I};$ \\
 (b) if $X\sim F$ has the moment sequence  $\{B_n\}_{n=0}^{\infty},$  then it has density function
 \begin{eqnarray}f(x)=\frac{1}{\pi}\frac{1}{\sqrt{x(4-x)}},\ x\in(0,4),
 \end{eqnarray} and moment generating function (mgf)
 \begin{eqnarray}M(t)=E[\exp(tX)]=\sum_{n=0}^{\infty}B_n\frac{t^n}{n!}=\exp(2t)I_0(2t),\ \forall\ t\in{\mathbb R},
 \end{eqnarray}
 where $I_{\alpha}(t)=\sum_{k=0}^{\infty}(t/2)^{2k+\alpha}/(k!\Gamma(k+\alpha+1)),\ t\in{\mathbb R},$
 is the modified Bessel function of the first kind;\\
(c) for each real $c>0,$ the power sequence
$\{B_n^c\}_{n=0}^{\infty}$ is
S-det;\\
(d) for each real $c>0,$ the measure $\mu_c$ corresponding to the
moment sequence $\{B_n^c\}_{n=0}^{\infty}$ has the Mellin transform
$${\cal
M}_c(s)=\int_0^{\infty}t^sd{\mu}_c(t)=\left(\frac{2^{2s}}{\sqrt{\pi}}\cdot\frac{\Gamma(s+1/2)}
{\Gamma(s+1)}\right)^c=\left(\frac{\Gamma(2s+1)}{(\Gamma(s+1))^2}\right)^c,\
s\ge 0.$$

Here are some interesting observations. It follows from (4) that
$Y=X/4$ has the arcsine density: $g(y)=\pi^{-1}(y(1-y))^{-1/2},\,
y\in(0,1)$. According to Example 3.5 in Berg and Dur\'an (2005), the
sequence $\{1/(n+1)\}_{n=0}^{\infty}$ belongs to ${\cal I}$, so
Theorem 1(a) is a consequence of Theorem 2(a) (by using Lemma 2
below). On the other hand, the next variant of Theorem 2 is used in
the proof of Theorem 3 below.

\noindent{\bf Theorem 2$^\prime$.} Let $\overline{B}_0=1,\
\overline{B}_n=$$2n\choose
n$/$2^n=\frac{1}{2^n}B_n=\frac{n+1}{2^n}C_n,\
n=1,2,\ldots.$ Then we have\\
(a) the real
 sequence $\{\overline{B}_n\}_{n=0}^{\infty}\in {\cal I};$ \\
  (b) if $Y\sim G$ has the moment sequence  $\{\overline{B}_n\}_{n=0}^{\infty},$
then it has density function
 $$g(y)=\frac{1}{\pi}\frac{1}{\sqrt{y(2-y)}},\ y\in(0,2),$$ and mgf
 $M(t)=E[\exp(tY)]=\exp(t)I_0(t),$ $\forall\ t\in{\mathbb R};$\\
(c) for each real $c>0,$ the power sequence
$\{\overline{B}_n^c\}_{n=0}^{\infty}$ is
S-det;\\
(d) for each real $c>0,$ the measure $\mu_c$ corresponding to the
moment sequence $\{\overline{B}_n^c\}_{n=0}^{\infty}$ has the Mellin
transform
$${\cal
M}_c(s)=\int_0^{\infty}t^sd{\mu}_c(t)=\left(\frac{2^{s}}{\sqrt{\pi}}\cdot\frac{\Gamma(s+1/2)}{\Gamma(s+1)}\right)^c
=\left(\frac{1}{2^s}\cdot \frac{\Gamma(2s+1)}{(\Gamma(s+1))^2}\right)^c,\,\
s\ge 0.$$

\noindent{\bf Theorem 3.} Let $D_0=1,\ D_n=(2n-1)!!, n=1,2,\ldots,$
where the double factorial number
$(2n-1)!!= (2n-1)(2n-3)\cdots 3\cdot 1.$ Then we have\\
(a) the double
factorial sequence $\{D_n\}_{n=0}^{\infty}\in {\cal I};$ \\
(b) for real $c>0,$ the power sequence $\{D_n^c\}_{n=0}^{\infty}$ is
S-det iff $c\le 2;$\\
(c) for each real $c\in (0,2],$ the measure $\mu_c$ corresponding to
the moment sequence $\{D_n^c\}_{n=0}^{\infty}$ has the Mellin
transform
$${\cal
M}_c(s)=\int_0^{\infty}t^sd{\mu}_c(t)=\left(\frac{2^s}{\sqrt{\pi}}\Gamma(s+1/2)\right)^c
=\left(\frac{1}{2^s}\cdot \frac{\Gamma(2s+1)}{\Gamma(s+1)}\right)^c,\ s\ge
0.$$

In the next result we consider the sequence
$\{(2n)!\}_{n=0}^{\infty}$ instead of the moment sequence
$\{n!\}_{n=0}^{\infty}$ in Lemma 5 below.

\noindent{\bf Theorem 4.}
(a) The moment sequence $\{(2n)!\}_{n=0}^{\infty}\in {\cal I}.$ \\
(b) For real $c>0,$  the power sequence
$\{((2n)!)^c\}_{n=0}^{\infty}$
is  S-det iff $c\le 1.$\\
(c) For each real $c\in(0,1],$ the measure $\mu_c$ corresponding to
the moment sequence $\{((2n)!)^c\}_{n=0}^{\infty}$ has the Mellin
transform ${\cal M}_c(s)=\int_0^{\infty}t^sd{\mu}_c(t)=
(\Gamma(2s+1))^c,$ $s\ge 0.$

For a fixed positive integer $k$ and the generalized Catalan numbers
\begin{eqnarray} C_{k,n}=\frac{1}{kn+1}{(k+1)n\choose n},\
n=0,1,2,\ldots \end{eqnarray} (or Fuss--Catalan numbers of order
$k$),  we have the following results.\\
\noindent{\bf Theorem 5.} Let $k\ge 2$ be a fixed positive integer.
 Then\\
(a) the sequence $\{C_{k,n}\}_{n=0}^{\infty}\in {\cal I};$ \\
(b) for each $c>0,$ the power sequence
$\{C_{k,n}^c\}_{n=0}^{\infty}$ is S-det;\\
(c) for each real $c>0,$ the measure $\mu_{k,c}$  corresponding to
the moment sequence $\{C_{k,n}^c\}_{n=0}^{\infty}$ has the Mellin
transform
$${\cal
M}_{k,c}(s)=\int_0^{\infty}t^sd{\mu}_{k,c}(t)
=\left(\frac{1}{ks+1}\cdot \frac{\Gamma((k+1)s+1)}{\Gamma(s+1)\Gamma(ks+1)}\right)^c,\
s\ge 0.$$
\noindent{\bf Theorem 6.}  Let $k\ge 2$ be a fixed positive
integer. Then\\
(a) the sequence  $\{{(k+1)n\choose n}\}_{n=0}^{\infty}\in {\cal I};$\\
(b) for each $c>0,$ the power sequence
$\{[{(k+1)n\choose n}]^c\}_{n=0}^{\infty}$ is S-det;\\
(c) for each real $c>0,$ the measure $\mu_{k,c}$  corresponding to
the moment sequence $\{[{(k+1)n\choose n}]^c\}_{n=0}^{\infty}$ has
the Mellin transform
$${\cal
M}_{k,c}(s)=\int_0^{\infty}t^sd{\mu}_{k,c}(t)
=\left(\frac{\Gamma((k+1)s+1)}{\Gamma(s+1)\Gamma(ks+1)}\right)^c,\
s\ge 0.$$ \noindent{\bf Theorem 7.} Let $k\ge 3$ be a fixed positive
integer.
Then \\
(a) $\{(kn)!\}_{n=0}^{\infty}$ is an {infinitely divisible} Stieltjes moment sequence;\\
(b) for real $c>0,$ the power sequence
$\{((kn)!)^c\}_{n=0}^{\infty}$ is S-det iff $kc\le 2.$
\bigskip\\
\noindent{\bf 3. Lemmas}

To prove the above main results, we need some notations and lemmas.
For a positive and infinitely differentiable function $f$ defined on
$(0,\infty),$ we say that $f$ is a completely monotone function if
$(-1)^nf^{(n)}(x)\ge 0$ for all $x\in (0,\infty)$ and $n\ge 1$,
denoted $f\in {\cal CM},$ and that $f$ is a Bernstein function,
denoted $f\in {\cal B},$ if its derivative $f^{\prime}\in {\cal
CM}.$ The following useful Lemmas 3--5 are due to Berg (2005, 2007),
Berg and Dur\'an (2004) as well as Berg and L\'opez (2015), in which
Lemma 3 extends a result of Bertoin and Yor (2001, Proposition 1).
\medskip\\
\noindent{\bf Lemma 1.} Let $\{m_n\}_{n=0}^{\infty}\in {\cal I}$ be
the moment sequence of a bounded nonnegative random variable $X.$
Then for each $c>0,$  the power sequence $\{m_n^c\}_{n=0}^{\infty}$ is S-det.\\
{\bf Proof.} Assume $X\sim F$ and $X\le M\ a.s.$ (almost surely),
where $M$ is a positive constant, and let $X_c\sim F_c$ be a
nonnegative random variable with moment sequence
$\{m_n^c\}_{n=0}^{\infty}.$ Then $m_n\le M^n$ and $m_n^c\le M^{cn}$
for all $n\ge 1$ and $c>0.$ The latter in turn implies that $X_c\le
M^c\ a.s.\ \forall\ c>0,$ because
$\lim_{n\to\infty}(E[X_c^n])^{1/n}=\lim_{n\to\infty}(m_n^c)^{1/n}\le
M^c.$ (See, e.g., Rudin 1987,  p.71, and note that
$E[X_c^n]=\int_0^1(F_c^{-1}(t))^ndt,$ where the quantile function
$F_c^{-1}(t)=\inf\{x: F_c(x)\ge t\},\  t\in(0,1).$) Therefore, $X_c$
is M-det on ${\mathbb R}_+$ by Hardy's or Cram\'er's criterion.
Another approach is to calculate the Carleman quantity for $X_c:$
$$C[X_c]\equiv\sum_{n=1}^{\infty}(E[X_c^n])^{-1/(2n)}=
\sum_{n=1}^{\infty}[m_n^c]^{-1/(2n)}\ge
\sum_{n=1}^{\infty}M^{-c/2}=\infty$$ and apply Carleman's
criterion (see, e.g., Lin 2017, Theorem 2). The proof is complete.
\medskip\\
 \noindent{\bf Lemma 2.} If the two Stieltjes moment sequences $\{s_n\}_{n=0}^{\infty}$ and
 $\{t_n\}_{n=0}^{\infty}$ are infinitely divisible, then so is the product sequence
$\{s_nt_n\}_{n=0}^{\infty}.$\\
{\bf Proof.} For each $c>0,$ let  $X_c$ and $Y_c$ be two independent
nonnegative random variables having the moment sequences
$\{s_n^c\}_{n=0}^{\infty}$ and $\{t_n^c\}_{n=0}^{\infty},$
respectively. Then the product $Z_c=X_cY_c$ has the moment sequence
$\{s_n^ct_n^c\}_{n=0}^{\infty},$ because
$E[Z_c^n]=E[X_c^n]E[Y_c^n]=s_n^ct_n^c\ \forall n$ due to the
independence of $X_c$ and $Y_c.$
 \medskip\\
 \noindent{\bf Lemma 3.} Let $\alpha,\,\beta>0$ and $f\in {\cal B}$
 with $f(\alpha)>0.$ Define the sequence $\{s_n\}_{n=0}^{\infty}$ by $s_0=1$ and
 $s_n=f(\alpha)f(\alpha+\beta)\cdots f(\alpha+(n-1)\beta),\ n=1,2,\ldots.$ Then  $\{s_n\}_{n=0}^{\infty}\in {\cal I},$
 and  the power sequence $\{s_n^c\}_{n=0}^{\infty}$ is S-det if $c\in (0,2].$
Moreover, for each real $c>0,$ let $X_c\sim F_c$ be a nonnegative
random variable having the moment sequence
$\{s_n^c\}_{n=0}^{\infty}.$ If $X_c$ is M-det, then the Mellin
transform of $F_c$ is of the form ${\cal
M}_c(z)\equiv\int_0^{\infty}t^zdF_c(t)=\exp[-c\psi(z)],$ Re$(z)\ge
0,$ where the exponent function
\begin{eqnarray*}\psi(z)=-z\log f(\alpha)+\int_0^{\infty}\left[(1-e^{-z\beta
x})-z(1-e^{-\beta x})\right]\frac{e^{-\alpha x}}{x(1-e^{-\beta
x})}d\kappa(x),
\end{eqnarray*}
and the measure $\kappa$ has the Laplace transform
\begin{eqnarray*}
\frac{f^{\prime}(s)}{f(s)}=\int_0^{\infty}e^{-sx}d\kappa(x),\ s> 0.
\end{eqnarray*}
 \noindent{\bf Lemma 4.} Assume that a Stieltjes moment sequence $\{m_n\}_{n=0}^{\infty}$
 is the product $m_n=s_nt_n$ of two Stieltjes moment sequences $\{s_n\}_{n=0}^{\infty}$ and
 $\{t_n\}_{n=0}^{\infty}.$ If $t_n>0$ for all $n$ and  $\{s_n\}_{n=0}^{\infty}$ is S-indet, then also
 $\{m_n\}_{n=0}^{\infty}$ is S-indet.
\medskip\\
\noindent{\bf Lemma 5.}
(a) The moment sequence $\{n!\}_{n=0}^{\infty}\in {\cal I}.$ \\
(b) For real $c>0,$  the power sequence $\{(n!)^c\}_{n=0}^{\infty}$
is  S-det iff $c\le 2.$\\
(c) For each real $c\in(0,2],$ the measure $\mu_c$ corresponding to
the moment sequence $\{(n!)^c\}_{n=0}^{\infty}$ has the Mellin
transform ${\cal M}_c(s)=\int_0^{\infty}t^sd{\mu}_c(t)=
(\Gamma(s+1))^c,$ $s\ge 0.$
\medskip\\
\noindent{\bf Lemma 6.} Let $X\sim F$ have the density $f$ defined
in (4), and let  $Z=X-2.$ Then \\
(a) the distribution of $Z$ is symmetric and $Z^2\stackrel{d}{=}X,$ where $\stackrel{d}{=}$ means `equal in distribution';\\
 (b) $Z$ has the mgf $M_Z(t)=E[\exp(tZ)]=I_0(2t),\ t\in{\mathbb R},$ where $I_0$ is
 the modified Bessel function of the first kind and is defined in (5);\\
 (c) the mgf of $X$ is equal to $M_X(t)=e^{2t}I_0(2t),\ t\in{\mathbb R}.$\\
{\bf Proof.} Clearly, $Z$ has a symmetric density function
\begin{eqnarray}
h(z)=\frac{1}{\pi}\frac{1}{\sqrt{4-z^2}},\ z\in(-2,2).
\end{eqnarray} We now
carry out the distribution of  $Z^2:$ for $x\in(0,4),$
\begin{eqnarray*}
\Pr(Z^2\le x)=\Pr(-\sqrt{x}<Z<\sqrt{x})=2\Pr(0\le Z<\sqrt{x})=\frac{2}{\pi}\int_0^{\sqrt{x}}\frac{1}{\sqrt{4-z^2}}dz.
\end{eqnarray*}
By changing variables $z=\sqrt{u},$ we have
$$\Pr(Z^2\le x)=\frac{1}{\pi}\int_0^{x}\frac{1}{\sqrt{u(4-u)}}du=F(x),\ \ x\in(0,4).$$
Therefore, $Z^2$ has the same distribution as $X.$ This proves part
(a). To prove part (b), we note that $E[Z^k]=0$ for odd $k$, and
$E[Z^{2n}]=E[X^n]=B_n={2n\choose n},\ n\ge 0.$ Hence,
$$M_Z(t)=E[\exp(tZ)]=\sum_{k=0}^{\infty}E[Z^k]\frac{t^k}{k!}=\sum_{n=0}^{\infty}E[Z^{2n}]
\frac{t^{2n}}{(2n)!}=\sum_{n=0}^{\infty}\frac{t^{2n}}{(n!)^2}=I_0(2t),\
\ t\in{\mathbb R}.$$ Alternatively, we can consider the random
variable $Z_*=Z/2$ which has the symmetric Beta distribution with
characteristic function
$$\phi_*(t)\equiv E[\exp(itZ_*)]=\sum_{n=0}^{\infty}\frac{(-1)^n}{(n!)^2}\left(\frac{t}{2}\right)^{2n},\ \ t\in{\mathbb R}$$
(see Kingman 1963, Theorem 1). This in turn implies that its mgf
$M_{Z_*}(t)=\phi_*(-it)=I_0(t),\ t\in{\mathbb R},$ and hence
$M_Z(t)=M_{Z_*}(2t)=I_0(2t),\ t\in{\mathbb R}.$ Finally, we have the
mgf of $X:$
$$M_X(t)=E[\exp(tX)]=E[\exp(t(Z+2))]=e^{2t}E[\exp(tZ)]=e^{2t}I_0(2t),\ \ t\in{\mathbb R}.$$ This completes the
proof.
\medskip\\
\indent For simplicity, we denote by $X\sim F \sim
\{m_n\}_{n=0}^{\infty}\sim {\cal M}$ the random variable $X$ having
distribution $F,$ moment sequence $\{m_n\}_{n=0}^{\infty}$ and
Mellin transform ${\cal M}.$
\medskip\\
\noindent{\bf Lemma 7.} Let $a$ and $c$ be two positive real
constants. Consider the random variables $X_1$ and $X_c$ satisfying
\begin{eqnarray*}
X_1\!&\sim&\! F_1 \sim \{m_n\}_{n=0}^{\infty}\sim {\cal M}_1,\ \ \ \ \ \ \ \ \ \ \
X_c\sim F_c \sim \{m_n^c\}_{n=0}^{\infty}\sim {\cal M}_c,\\
X_1/a\!&\sim&\! F_{1,a} \sim \{m_n/a^n\}_{n=0}^{\infty}\sim {\cal M}_{1,a},\ \
X_c/a^c\sim F_{c,a} \sim \{(m_n/a^n)^c\}_{n=0}^{\infty}\sim {\cal M}_{c,a},
\end{eqnarray*}
and ${\cal M}_c(s)={\cal M}_1^c(s),\ s\ge 0.$ Then  the Mellin
transform
${\cal M}_{c,a}(s)={\cal M}_{1,a}^c(s),\ s\ge 0.$\\
{\bf Proof.} For $s\ge 0,$ we have
\begin{eqnarray*}
{\cal M}_{c,a}(s)&=&E[(X_c/a^c)^s]=\frac{1}{a^{cs}}E[X_c^s]=\frac{1}{a^{cs}}{\cal M}_c(s)=\frac{1}{a^{cs}}{\cal M}_1^c(s)\\
&=&\left(\frac{1}{a^{s}}E[X_1^s]\right)^c=\left(E[(X_1/a)^s]\right)^c={\cal M}_{1,a}^c(s).
\end{eqnarray*}
The proof is complete.
\bigskip\\
{\bf 4. Proofs of main results}\medskip\\
{\bf Proof of Theorem 1.} Rewrite first  \begin{eqnarray*}
C_n=\frac{(2n)!}{n!(n+1)!}&=&\frac{2n}{n}\cdot\frac{2n-1}{n+1}\cdot\frac{2n-2}{n-1}
\cdot\frac{2n-3}{n}\cdots\cdot\frac{5}{4}\cdot\frac{4}{2}\cdot\frac{3}{3}\cdot\frac{2}{1}\cdot\frac{1}{2}\\
&=&2^n\cdot\frac{2n-1}{n+1}\cdot\frac{2n-3}{n}\cdots\cdot\frac{5}{4}\cdot\frac{3}{3}\cdot\frac{1}{2}\\
&=&2^n\left(2-\frac{3}{n+1}\right)\left(2-\frac{3}{n}\right)\cdots\cdot
\left(2-\frac{3}{4}\right)\left(2-\frac{3}{3}\right)\left(2-\frac{3}{2}\right).
\end{eqnarray*}
Then set the functions $h(x)=2(2-3/(x+1)),\ x>1/2,$ and
$f(x)=h(x+1/2),\ x> 0.$ Next, take $\alpha=1/2$ and $\beta=1$. We
have $f(\alpha)=h(1)=1>0, f\in {\cal B},$ and the Catalan number
$C_n=f(\alpha)f(\alpha+\beta)\cdots
f(\alpha+(n-1)\beta)=\Pi_{k=1}^nh(k)$ for all $n\ge 1.$ Lemma 3
applies and hence the Catalan sequence $\{C_n\}_{n=0}^{\infty}\in
{\cal I}.$ This proves part (a). Part (b) follows from Lemma 1 and
the fact (3) immediately. To prove part (c), we write, by Lemma 3,
$${\cal
M}_c(s)\equiv\int_0^{\infty}t^sd{\mu}_c(t)=\left[\int_0^{\infty}t^sd{\mu}_1(t)\right]^c={\cal
M}_1^c(s),\ s\ge 0.$$ Therefore, it remains to carry out ${\cal M}_1$ and we have,  by (3), that, for $s\ge 0,$
\begin{eqnarray*}
{\cal M}_1(s)&=&\int_0^{\infty}t^sd{\mu}_1(t)=\int_0^4t^s{\left(\frac{1}{2\pi}\sqrt{\frac{4-t}{t}}\right)}dt\\
&=&\frac{2^{2s+1}}{\pi}\int_0^1x^{s-1/2}(1-x)^{1/2}dx=\frac{2^{2s+1}}{\pi}\cdot\frac{\Gamma(s+1/2)\Gamma(3/2)}{\Gamma(s+2)}\\
&=&\frac{2^{2s}}{\sqrt{\pi}}\cdot\frac{\Gamma(s+1/2)}{\Gamma(s+2)}=\frac{1}{s+1}\cdot \frac{\Gamma(2s+1)}{(\Gamma(s+1))^2}.
\end{eqnarray*}
The last equality is due to the duplication formula for the gamma
function: $$\Gamma(x)\Gamma(x+1/2)=2^{1-2x}\sqrt{\pi}\,\Gamma(2x),\,
x>0$$ (see, e.g.,  Gradshteyn and Ryzhik 2014, Section 8.335). The
proof is complete.
\medskip\\
{\bf Proof of Theorem 2.} Rewrite first \begin{eqnarray*}
B_n=\frac{(2n)!}{(n!)^2}&=&\frac{2n}{n}\cdot\frac{2n-1}{n}\cdot\frac{2n-2}{n-1}
\cdot\frac{2n-3}{n-1}\cdots\cdot\frac{5}{3}\cdot\frac{4}{2}\cdot\frac{3}{2}\cdot\frac{2}{1}\cdot\frac{1}{1}\\
&=&2^n\cdot\frac{2n-1}{n}\cdot\frac{2n-3}{n-1}\cdots\cdot\frac{5}{3}\cdot\frac{3}{2}\cdot\frac{1}{1}\\
&=&2^n\left(2-\frac{1}{n}\right)\left(2-\frac{1}{n-1}\right)\cdots\cdot
\left(2-\frac{1}{3}\right)\left(2-\frac{1}{2}\right)\left(2-\frac{1}{1}\right).
\end{eqnarray*}
Set the functions $h(x)=2(2-1/x),\ x>1/2,$ and $f(x)=h(x+1/2),\
x>0.$ It is seen that $f\in {\cal B}.$ Next, take $\alpha=1/2$ and
$\beta=1,$  then we have $f(\alpha)=2>0$ and
$B_n=f(\alpha)f(\alpha+\beta)\cdots
f(\alpha+(n-1)\beta)=\Pi_{k=1}^nh(k),\ n\ge 1.$ Lemma 3 applies and
this proves part (a). To prove part (b), assume that $U$ is a
uniform random variable with values in $(0,1)$ and is independent of
$X.$ Then the product $XU$ has the moment sequence
$\{C_n\}_{n=0}^{\infty}$  and hence has the density $f_C$ defined in
(3) because the Catalan sequence is S-det. Therefore, we have
$$\Pr(XU\le y)=\int_0^{y}f_C(x)dx,\ y\ge 0,$$
or, equivalently,
\begin{eqnarray}\int_0^1\Pr(X\le
y/u)du=\int_0^y\frac{1}{2\pi}\sqrt{\frac{4-x}{x}}dx,\ y\in (0,4).
\end{eqnarray}By changing variables $y/u=x,$ we rewrite (8) as
\begin{eqnarray}
\int_y^{\infty}F(x)\frac{y}{x^2}dx=\int_0^y\frac{1}{2\pi}\sqrt{\frac{4-x}{x}}dx,\ y\in (0,4).
\end{eqnarray}
Solving Eq.\,(9) by differentiating both sides twice with respect to
$y,$ we prove the density function $f$ of $X$ to be the one given in
(4). The identity in (5) is proved in Lemma 6.

Part (c) follows from parts (a) and (b) and Lemma 1 immediately.
Finally, we carry out the Mellin transform of $X$ (or $F$). For each
$s\ge 0,$ $E[(XU)^s]=E[X^s]E[U^s]=E[X^s]/(s+1)$ happens to be the
Mellin transform of $\mu_1$ corresponding to the Catalan sequence,
namely, $\frac{1}{s+1}\frac{\Gamma(2s+1)}{(\Gamma(s+1))^2},$ from
which the required $E[X^s]$ follows.  This completes the proof.
\medskip\\
{\bf Proof of Theorem 2$^\prime$.} (Method I) Note that if
$\{s_n\}_{n=0}^{\infty}\in{\cal I},$ then so is
$\{a^ns_n\}_{n=0}^{\infty}$ for each constant $a>0.$ Therefore,
Theorem 2$^\prime$ follows from Theorem 2, and vice versa.

(Method II) Define the functions $h(x)=2-1/x,\ x>1/2,$ and
$f(x)=h(x+1/2),\ x>0.$ Take $\alpha=1/2$ and $\beta=1,$ then
 $f\in {\cal B}$ and $f(\alpha)=1>0.$ The remaining proof is similar to that of
Theorem 2. Note also that $Y$ has the same distribution as $X/2,$
and hence the required results follow from Theorem 2 and Lemma 7 as
well. The proof is complete.
\medskip\\
{\bf Proof of Theorem 3.} Consider the functions $h(x)=2x-1,\
x>1/2,$ and $f(x)=h(x+1/2),\ x>0.$ It is seen that $f\in {\cal B}.$
Next, take $\alpha=1/2$ and $\beta=1,$  then we have $f(\alpha)=1>0$
and $D_n=f(\alpha)f(\alpha+\beta)\cdots
f(\alpha+(n-1)\beta)=\Pi_{k=1}^nh(k),\ n\ge 1.$ Lemma 3 applies and
this proves both (a) and the sufficient part of (b). To prove the
necessary part of (b), we write first $D_n=n!\cdot (2n)!/[2^n
(n!)^2]=n!\cdot \overline{B}_n\equiv s_n\cdot t_n$ for all $n\ge 1$
and $D_n^c=s_n^c\cdot t_n^c$ for $c>0$ (see Theorem 2$^\prime$).
Namely, the Stieltjes moment sequence $\{D_n^c\}_{n=0}^{\infty}$ is
the product of two Stieltjes moment sequences
$\{s_n^c\}_{n=0}^{\infty}$ and $\{t_n^c\}_{n=0}^{\infty}.$ If $c>2,$
$\{s_n^c\}_{n=0}^{\infty}$ is S-indet (Lemma 5), and hence
$\{D_n^c\}_{n=0}^{\infty}$ is also S-indet by Lemma 4.

To prove (c), let us recall that $\{D_n\}_{n=0}^{\infty}$  is in
fact the moment sequence of the chi-square random variable
$\chi_1^2$ (the square of the standard normal random variable) with
density function $f(x)=\frac{1}{\sqrt{2\pi}}x^{-1/2}\exp(-x/2),\
x>0,$ and Mellin transform
$${\cal M}_1(s)=\int_0^{\infty}t^sf(t)dt=\frac{1}{\sqrt{2\pi}}\int_0^{\infty}t^{s-1/2}e^{-t/2}dt
=\frac{2^s}{\sqrt{\pi}}\Gamma(s+1/2),\ s\ge 0.$$ The proof is complete.
\medskip\\
{\bf Proof of Theorem 4.} Let the random variables ${\cal E}$ and
$X$ have the moment sequences $\{n!\}_{n=0}^{\infty}$ and
$\{B_n\}_{n=0}^{\infty}=\{{2n\choose n}\}_{n=0}^{\infty},$
respectively.  Then $Z={\cal E}^2$ has the moment sequence
$s_n\equiv E[Z^n]=(2n)!,\ n\ge 0.$ Next, write $s_n={2n\choose
n}(n!)^2=B_n(n!)^2$ and $s_n^c=B_n^c(n!)^{2c}$ for each $c>0.$
Therefore, $\{s_n\}_{n=0}^{\infty}\in {\cal I}$ by Theorem 2(a) and
Lemma 5(a). This proves part (a). Moreover, if $c>1,$ the sequence
$\{(n!)^{2c}\}$ is S-indet due to Lemma 5(b), so is
$\{s_n^c\}_{n=0}^{\infty}$  by Lemma 4. On the other hand, if
$c\in(0,1],$ the Carleman quantity
$\sum_{n=0}^{\infty}[s_n^c]^{-1/(2n)}=\infty$ by Stirling formula
$\Gamma(x+1)=x\Gamma(x)\cong \sqrt{2\pi}x^{x+1/2}e^{-x}$ as
$x\to\infty,$ and hence $\{s_n^c\}_{n=0}^{\infty}$ is S-det. This
completes the proof of part (b). To prove part (c), let $c\in(0,1],$
and let $X_c,\ {\cal E}_{2c}$ be two independent nonnegative random
variables having moment sequences $\{B_n^c\}_{n=0}^{\infty},\
\{(n!)^{2c}\}_{n=0}^{\infty},$ respectively. Finally, define
$Z_c=X_c {\cal E}_{2c},$ which has the moment sequence
$\{s_n^c\}_{n=0}^{\infty}.$ Then we have, by Theorem 2(d) and Lemma
5(c), that
$$E[Z_c^s]=E[X_c^s]E[{\cal E}_{2c}^s]=\left(\frac{\Gamma(2s+1)}{(\Gamma(s+1))^2}\right)^c
(\Gamma(s+1))^{2c}=(\Gamma(2s+1))^c,\ s\ge 0.$$ The proof is complete.
\medskip\\
{\bf Proof of Theorem 5.} Mimicking the proof of Theorem 1, we
define for the sequence $\{C_{k,n}\}_{n=0}^{\infty}$ the functions:
$$h_{\ell}(x)=(k+1)^{1/k}\left[(1+1/k)-\frac{2-(\ell-2)/k}{kx-(\ell-2)}\right],\ x>\ell/(k+1),\
\ell=1,2,\ldots,k.$$ Next, define the functions:
$f_{\ell}(x)=h_{\ell}(x+\ell/(k+1)),\ x>0,\ \ell=1,2,\ldots,k.$ Note
that all the functions $f_{\ell}\in{\cal B}.$ Taking
$\alpha_{\ell}=1-\ell/(k+1),\ \beta_{\ell}=1,\ \ell=1,2,\ldots,k,$
we have
$f_{\ell}(\alpha_{\ell})=h_{\ell}(1)=(k+1)^{1/k}\frac{k-\ell+1}{k-\ell+2}>0,\
\ell=1,2,\ldots,k,$ and
$$C_{k,n}=\Pi_{\ell=1}^k\left(\Pi_{i=1}^nf_{\ell}(\alpha_{\ell}+(i-1)\beta_{\ell})\right)=\Pi_{\ell=1}^k\left(\Pi_{i=1}^nh_{\ell}(i)\right).$$
By Lemma 3, each sequence
$\{\Pi_{i=1}^nh_{\ell}(i)\}_{n=0}^{\infty}\in{\cal I}$ (with the
first term equal to $1$), so is the product sequence
$\{C_{k,n}\}_{n=0}^{\infty}$ due to Lemma 2. This proves part (a).
To prove part (b), let us
 recall that the  measure corresponding to
the moment sequence $\{C_{k,n}\}_{n=0}^{\infty}$ has a bounded
support $(0, (k+1)^{k+1}/k^k]$ and hence Lemma 1 applies (see, e.g.,
Banica et al.\,2011, Theorem 2.1).

Finally,  to prove part (c), we recall first that the Mellin
transform of the measure $\mu_{k,1}$ is
$${\cal
M}_{k,1}(s)=\int_0^{\infty}t^sd{\mu}_{k,1}(t) =\frac{1}{ks+1}\cdot
\frac{\Gamma((k+1)s+1)}{\Gamma(s+1)\Gamma(ks+1)},\ s\ge 0$$ (see
Penson and $\dot{\hbox{Z}}$yczkowski 2011, p.\,2). Then Lemma 3
completes the proof.
\medskip\\
{\bf Proof of Theorem 6.} As before, mimicking the proof of Theorem
2, we define the functions:
$$h_{\ell}^*(x)=(k+1)^{1/k}\left[(1+1/k)-\frac{(k-\ell+1)/k}{kx-(\ell-1)}\right],\ x>\ell/(k+1),\
\ell=1,2,\ldots,k,$$ and $f_{\ell}^*(x)=h_{\ell}^*(x+\ell/(k+1)),\
x>0,\ \ell=1,2,\ldots,k.$ All the functions $f_{\ell}^*\in{\cal B}.$
Taking $\alpha_{\ell}=1-\ell/(k+1),\ \beta_{\ell}=1,\
\ell=1,2,\ldots,k,$ we have
$f_{\ell}^*(\alpha_{\ell})=h_{\ell}^*(1)=(k+1)^{1/k}>0,\
\ell=1,2,\ldots,k,$ and
$$(kn+1)C_{k,n}={(k+1)n\choose n}=\Pi_{\ell=1}^k\left(\Pi_{i=1}^n
f_{\ell}^*(\alpha_{\ell}+(i-1)\beta_{\ell})\right)=\Pi_{\ell=1}^k\left(\Pi_{i=1}^nh_{\ell}^*(i)\right).$$
Again, Lemmas 2 and 3 apply, and the sequence $\{{(k+1)n\choose
n}\}_{n=0}^{\infty}\in{\cal I}.$ This proves part (a).

To prove part (b), we can use Cram\'er's or Carleman's criterion
again. For the rest of the proof, let $Y$ be a nonnegative random
variable having the moment sequence
$\{(kn+1)C_{k,n}\}_{n=0}^{\infty},$ and let $U$ be a uniform random
variable on $(0,1)$ and independent of $Y.$ Then the product
$X=YU^k$ has the moment  sequence $\{C_{k,n}\}_{n=0}^{\infty}$ and
Mellin transform $E[X^s]=E[Y^sU^{ks}]=\frac{1}{ks+1}E[Y^s],\,s\ge
0.$ This implies by Theorem 5(c) that
$$E[Y^s]=\frac{\Gamma((k+1)s+1)}{\Gamma(s+1)\Gamma(ks+1)},\ \ s\ge 0,$$
and finally, Lemma 3 completes the proof of part (c).
\medskip\\
{\bf Proof of Theorem 7.} Rewrite
$$(kn)!={{kn}\choose n}{{(k-1)n}\choose n}\cdots{{2n}\choose n}(n!)^k\equiv A_k(n)A_{k-1}(n)\cdots A_2(n)A_1(n).$$
Recall that each $\{A_i(n)\}_{n=0}^{\infty}\in{\cal I}$ by Lemma 5
and Theorems 2 and 6, where $i=1,2,\ldots k,$ so is the product
sequence $\{(kn)!\}_{n=0}^{\infty}$ due to Lemma 2. This proves (a).
To prove the sufficiency part of (b), we can use Carleman's
criterion and Stirling formula, while the necessity part of (b)
follows from the result of (a) as well as Lemmas 4 and 5(b). The
proof is complete.
\bigskip\\
{\bf 5. Remarks}

\noindent{\bf Remark 1.} The property in Lemma 6(a) happens to be a
characterization of the density (4). In other words, if the
distribution of $Z$ is symmetric and $Z^2\stackrel{d}{=}Z+2,$ then
$Z$ has the density (7). To see this, we have from the above
equality that $Z\ge -2$ and hence $-2\le Z\le 2$ due to the
symmetric condition. This further implies that $Z$ is M-det. By the
equality in distribution, we can carry out all the moments of $Z$
and claim that $Z$ has the density (7).

\noindent{\bf Remark 2.} Recall that the mgf of Beta distribution
$B_{\alpha,\beta}$ is
$$M_{\alpha,\beta}(t)=1+\sum_{n=1}^{\infty}\left(\Pi_{r=0}^{n-1}
\frac{\alpha+r}{\alpha+\beta+r}\right)\frac{t^n}{n!},\ \
t\in{\mathbb R}.$$ It can be shown that the mgf of the Catalan
distribution $F_C$ with density $f_C$ in (3) is equal to
$M_C(t)=M_{\frac{1}{2},\frac{3}{2}}(4t),\ \forall t.$ On the other
hand, with the help of Theorem 2(b) and the definition of mgf,  we
have the relation: $\frac{d}{dt}\{tM_C(t)\}=\exp(2t)I_0(2t),\
t\in{\mathbb R},$ from which it follows that
$M_C(t)=\exp(2t)[I_0(2t)-I_1(2t)],\ t\in{\mathbb R}.$

\noindent{\bf Remark 3.} The Fuss--Catalan density (of order 2)
corresponding to the moment sequence $\{C_{2,n}\}_{n=0}^{\infty}$ in
(6) is of the form
\begin{eqnarray*}
f_2(x)
=\frac{3\left(1+\sqrt{1-4x/27}\right)^{2/3}-2^{2/3}x^{1/3}}
{2^{4/3}3^{1/2}\pi x^{2/3}\left(1+\sqrt{1-4x/27}\right)^{1/3}},\ \ x\in(0,27/4],
\end{eqnarray*}
 due to Penson
and Solomon (2001) (see also M{\l}otkowski et al.\,2013, Theorem
4.3). For the complicated general case, see Liu et al.\,(2011) or
Penson and $\dot{\hbox{Z}}$yczkowski (2011).

 \noindent{\bf Remark 4.} Recall that the powers of the
standard exponential random variable ${\cal E}$ have the moment
(in)determinacy property: for each real $c>0,$ ${\cal E}^c$ is M-det
iff $c\le 2$ (see, e.g., Targhetta 1990). This together with Lemma 5
implies that for each $c>0,$ the $c$-th power  ${\cal E}^c$ of
${\cal E}$ is M-det iff the $c$-th power $\{(n!)^c\}_{n=0}^{\infty}$
of $\{n!\}_{n=0}^{\infty}$ is S-det.

\noindent{\bf Remark 5.} It is known that the mgf of the chi-square
random variable $\chi_1^2$ (with mean 1) is $M(t)=1/\sqrt{1-2t},\
t<\frac{1}{2}.$ The powers of $\chi_1^2$ share the same M-(in)det
property with those of the standard exponential random variable
${\cal E},$ namely, $(\chi_1^2)^c$ is M-det iff $c\le 2$ (see, e.g.,
Berg 1988). This together with Theorem 3 implies that for each
$c>0,$ the $c$-th power $(\chi_1^2)^c$ of $\chi_1^2$    is M-det iff
the $c$-th power $\{D_n^c\}_{n=0}^{\infty}$ of the chi-square moment
sequence $\{D_n\}_{n=0}^{\infty}$
  is S-det. An interesting by-product in
the proof of Theorem 3 is the factorization of $\chi_1^2.$ More
precisely, $\chi_1^2$ has the same distribution with the product
${\cal E}\cdot Y,$ where the random variable $Y$ is defined in
Theorem 2$^{\prime}$ and is independent of ${\cal E}.$

\noindent{\bf Remark 6.} Recall also that any positive real power of
a bounded nonnegative random variable is M-det. In view of these
existing results including Theorems 1 and 3 as well as Lemma 5, we
would like to pose the following.

\noindent{\bf Conjecture 1.} Let $X$ be a nonnegative random
variable with finite moments $m_n=E[X^n]$ of all orders $n\ge 1.$
Moreover, let $\{m_n\}_{n=0}^{\infty}\in {\cal I}$ and for each real
$c>0,$ denote by $X_c$ a nonnegative random variable having the
moment sequence $\{m_n^c\}_{n=0}^{\infty}.$ Then the $c$-th power
$X^c$  is M-det on ${\mathbb R}_+$ iff  the $c$-th power sequence
$\{m_n^c\}_{n=0}^{\infty}$ is S-det, equivalently, $X_c$ is M-det on
${\mathbb R}_+.$

\noindent{\bf Remark 7.} Note that for each real $a>0,$
$\{\Gamma(an+1)\}_{n=0}^{\infty}$ is a Stieltjes moment sequence
because $E[({\cal E}^a)^n]=\Gamma(an+1)$ for each $n\ge 0.$ In view
of Lemma 5 and Theorems 4 and 7, we also pose the following.

\noindent{\bf Conjecture 2.} Let $a>0$ be a real constant and
$m_n=\Gamma(an+1),\ n\ge 0.$
Then\\
(a) the moment sequence
$\{m_n\}_{n=0}^{\infty}\in{\cal I};$\\
(b) for real $c>0,$ the power sequence $\{m_n^c\}_{n=0}^{\infty}$ is
S-det iff $ac\le 2;$\\
(c) for each real $c\in (0,2/a],$ the measure $\mu_c$ corresponding
to the moment sequence $\{m_n^c\}_{n=0}^{\infty}$ has the Mellin
transform ${\cal M}_c(s)=\int_0^{\infty}t^sd{\mu}_c(t)=
(\Gamma(as+1))^{c},$ $s\ge 0.$
\smallskip\\
\noindent{\bf Acknowledgments.} The author would like to thank the
 two Referees for helpful comments and constructive suggestions
which improved the presentation of the paper. {\it In particular,
one Referee points out that Conjecture 2 and its generalization have
been established in the recent article Berg (2018).}
\bigskip\\
\centerline{\bf References}
\begin{description}

\item {Akhiezer, N.\,I. (1965)}. {\it The Classical Problem of
Moments and Some Related Questions of Analysis}. Oliver $\&$
Boyd, Edinburgh. [Original Russian edn: Nauka,  Moscow, 1961.]

\item Amdeberhan, T., Moll, V. and Vignat, C. (2013). A probabilistic
interpretation of a sequence related to Narayana polynomials.
{\it Online J. Anal. Comb.}, {\bf 8}, 1--25.

\item Banica, T., Belinschi, S.\,T., Capitaine, M. and Collins, B. (2011).
Free Bessel laws. {\it Canad. J. Math.}, {\bf 63},  3--37.

\item Berg, C.
(1988). The cube of a normal distribution is indeterminate. {\it
Ann. Probab.} {\bf 16}, 910--913.

\item Berg, C. (2005). On powers of Stieltjes moment sequences,
I. {\it J. Theor. Probab.}, {\bf 18}, 871--889.

\item Berg, C. (2007). On powers of Stieltjes moment sequences,
II. {\it J. Comput. Appl. Math.}, {\bf 199}, 23--38.

\item Berg, C. (2018). A two-parameter extension of the Urbanik semigroup. arXiv:1802.00993.

\item Berg, C. and Dur\'an, A.\,J. (2004). A transformation from Hausdorff to Stieltjes moment sequences.
{\it Ark. Mat.}, {\bf 42}, 239--257.

\item Berg, C. and Dur\'an, A.\,J. (2005). Some transformations of Hausdorff moment sequences and harmonic numbers.
{\it Canad. J. Math.}, {\bf 57}, 941--960.

\item Berg, C. and L\'opez, J.\,L. (2015). Asymptotic behaviour of the Urbanik semigroup. {\it J. Approx. Theory}, {\bf 195}, 109--121.

\item Bertoin, J. and Yor, M. (2001). On subordinators, self-similar Markov processes and some factorizations of the exponential variables.
{\it Elect. Comm.  Probab.}, {\bf 6}, 95--106.

\item Binet, J. (1839). R\'eflexions sur le probl$\grave{\hbox{e}}$me de d\'eterminer le nombre de mani$\grave{\hbox{e}}$res dont une figure
rectiligne peut $\hat{\hbox{e}}$tre partag\'ee en triangles au
moyen de ses diagonals. {\it J. Math. Pures Appl.}, {\bf 4},
79--91.

\item Catalan, E. (1838). Note sur une \'equation aux diff\'erences. {\it J. Math. Pures Appl.}, {\bf 3},
508--516.

\item Gradshteyn, I.\,S. and Ryzhik, I.\,M. (2014). {\it Table of Integrals, Series, and Products}. 8th edn.  Academic Press, New York.

\item Kingman, J.\,F.\,C. (1963). Random walks with spherical symmetry. {\it Acta Math.}, {\bf 109}, 11--53.

\item Koshy, T. (2009). {\it Catalan Numbers with Applications.} Oxford University Press, Oxford.

\item Larcombe, P. (1999). The 18th century Chinese discovery of the Catalan numbers. {\it Mathematical
Spectrum}, {\bf 32}, 5--7.

\item  Liang, H., Mu, L. and Wang, Y. (2016). Catalan-like
numbers and Stieltjes moment sequences. {\it Discrete Math.},
{\bf 339}, 484--488.

\item Lin, G.\,D. (2017). Recent developments on the moment problem. {\it Journal of
 Statistical Distributions and  Applications}, {\bf 4}:5.

 \item Liu, D.-Z., Song, C. and Wang, Z.-D. (2011). On explicit probability densities associated with Fuss--Catalan numbers.
 {\it Proc. Amer. Math. Soc.}, {\bf 139}, 3735--3738.

\item Luo, J.\,J. (1988). Antu Ming, the first discoverer of the Catalan numbers (in Chinese).
{\it Neimenggu Daxue Xuebao},  {\bf 19},  239--245.

\item Luo, J.\,J. (2013). Ming Antu and his power series expansions. In:
{\it Seki, Founder of Modern Mathematics in Japan: A
 Commemoration on His Tercentenary}. E. Knobloch et al. (eds).
 Springer Proc. Math. Stat., {\bf 39}, 299--310, Springer,
 Tokyo.

\item M{\l}otkowski, W., Penson, K.\,A. and   $\dot{\hbox{Z}}$yczkowski, K. (2013).
Densities of the Raney distributions. {\it Documenta Math.},
{\bf 18}, 1573--1596.

 \item Penson,  K.\,A. and Solomon, A.\,I. (2001). Coherent states from combinatorial sequences.
 In: {\it Quantum Theory and Symmetries} (Krakow, 2001), pp.
 527--530.  E. Kapuscik and A. Horzela (eds), World Sci. Publ.,
 River Edge, NJ, 2002.

 \item Penson,  K.\,A. and $\dot{\hbox{Z}}$yczkowski, K. (2011). Product of Ginibre matrices: Fuss--Catalan and Raney distributions.
 {\it Phys. Rev.} E, {\bf 83}, 1--9.

 \item Roman, S. (2015). {\it An Introduction to Catalan Numbers}. Birkh\"auser, New York.

 \item Rudin, W. (1987). {\it Real and Complex Analysis.} 3rd edn. McGraw-Hill, New York.

\item {Shohat, J.\,A. and Tamarkin, J.\,D. (1943)}.
{\it The Problem of Moments.}  Amer. Math. Soc., New York.

\item Stanley, R.\,P. (1999). {\it Enumerative Combinatorics}, Vol. II. Cambridge University Press, Cambridge.

\item Stanley, R.\,P. (2015).  {\it Catalan Numbers}. Cambridge University Press, Cambridge.

\item {Stieltjes, T.\,J.} (1894/1895). Recherches sur les fractions
continues. {\it Ann. Fac. Sci. Univ. Toulouse Math.} {\bf 8(J)},
1--122 (1894); {\bf 9(A)}, 1--47 (1895). Also in: {Stieltjes,
T.J. Oeuvres Completes}. Noordhoff, Gr\"oningen {\bf 2},
402--566 (1918)

\item Targhetta,  M.\,L. (1990). On a family of indeterminate distributions. {\it J. Math. Anal. Appl.}, {\bf 147}, 477--479.

\item Tyan, S.-G. (1975). The structure of bivariate distribution functions and their relation to Markov processes.
Ph.D.\,Thesis, Princeton University, USA.

\end{description}
\end{document}